\documentclass{amsart}
\usepackage{url}
\usepackage{mathtools}

\usepackage[citecolor=blue,urlcolor=magenta,linkcolor=true, hypertexnames=false]{hyperref}

\theoremstyle{definition}

\theoremstyle{remark}

\numberwithin{equation}{section}

\DeclarePairedDelimiter\abs{\lvert}{\rvert}

\DeclarePairedDelimiter\ceil{\lceil}{\rceil}

\usepackage{xcolor}

\begin{document}

\title[Improved algorithms for left factorial residues]{Improved algorithms
for left factorial residues}

\author[Andreji\'c]{Vladica Andreji\'c$^\dagger$}
\thanks{$^\dagger$ Faculty of Mathematics, University of Belgrade, Serbia; \texttt{andrew@matf.bg.ac.rs}}

\author[Bostan]{Alin Bostan$^\star$}
\thanks{$^\star$ Inria, Univ. Paris-Saclay, France; \texttt{alin.bostan@inria.fr}}           

\author[Tatarevic]{Milos Tatarevic$^\ddagger$}
\thanks{$^\ddagger$ Alameda, CA 94501; \texttt{milos.tatarevic@gmail.com}}


\date{}

\dedicatory{}

\begin{abstract} We present improved algorithms for computing the left
factorial residues $!p=0!+1!+\dots+(p-1)! \!\mod p$. We use these algorithms
for the calculation of the residues $!p\!\mod p$, for all primes $p$~up to
$2^{40}$. Our results confirm that Kurepa's \emph{left factorial conjecture}
is still an open problem, as they show that there are no odd primes~$p<2^{40}$
such that $p$ divides $!p$. Additionally, we confirm that there are no
\emph{socialist primes}~$p$ with $5<p<2^{40}$. \end{abstract}

\maketitle

\section{Introduction}

In 1971, Kurepa~\cite{K1} introduced the left factorial function $!n$ defined, for any positive integer~$n$,  as the sum of factorials
$!n=0!+1!+\dots+(n-1)!$. Kurepa conjectured that the greatest common divisor
of $!n$ and $n!$ is equal to $2$ for all integers $n>2$. Equivalently, the
conjecture claims that there are no odd primes $p$ such that $p$ divides $!p$.
This problem has been studied extensively and was called \emph{Kurepa's
conjecture} by the subsequent authors. For a historical background, the reader
can consult~\cite{AT1}. The conjecture is also listed in Richard Guy's
classical book~\cite[Section~B44]{Gu}; as of 2020, it is still an
open problem.

In the past, there were several attempts to disprove the conjecture by finding
a counterexample. In the most recent search~\cite{AT1}, no counterexample was
found for any $p < 2^{34}$. All such searches are based on calculations of
residues $r_p = \,\,!p\!\mod p$ for primes $p$. In all previous attempts, the
time complexity of algorithms was $O(p)$ for a single $p$ and $O(n^2 /\log n)$
for all $p<n$. We now show that the computational complexity can be
significantly improved and we extend the search range up to $2^{40}$.

These improvements are based on the simple observation that $n!$ and $!n$
can be represented altogether in a matrix factorial form as
\begin{equation*}
M_n := C_1 C_2 \cdots C_n =\begin{pmatrix}n! & !n\\ 0 & 1\end{pmatrix},
\qquad \text{where} \quad
C_k=\begin{pmatrix}k & 1\\ 0 & 1\end{pmatrix}.
\end{equation*}

Applying \cite[Theorem 8]{BGS} on the matrix factorial $M_p$, yields an
improved algorithm for computing a single remainder $r_p$ in time
$O(p^{0.5+\varepsilon})$. In practice, due to the overhead of the 
polynomial multiplication based on Fast Fourier Transform, this method does not significantly outperform the
one given in \cite{AT1} for $p$ around $2^{34}$. However, the improvement is
notable for larger values of $p$. We used this method to verify the individual
values $r_p$ obtained by the algorithm we will describe next.

The main method we used in our search is based on the work presented in
\cite{CGH}. The algorithm was originally designed for computing \emph{Wilson
primes}, and is easy to adapt to matrix factorials. As a consequence, the
remainders $r_p$ for $2\leq p\leq n$ can be computed altogether in time
$O(n\log^{3+\varepsilon}n)$. However, for large $n$, due to the limited
computing resources, we need to run the search on smaller intervals. The time
complexity of the method we used to compute the remainders $r_p$ for $m\leq
p\leq n$ is 
\begin{equation*} 
O(m\log^{3+\varepsilon}m + (n-m)\log^{3+\varepsilon}n). 
\end{equation*}

Let us denote the terms we will use in the following section, where some definitions
are similar to those that appear in \cite[Theorem 2]{CGH}. 
Let $h = \ceil{\log_2 (n-m)}$. For each $0\leq i \leq h$ and $0\leq j\leq 2^i$ we set
\begin{equation*}
S_{i,j}=\left\{k\in\mathbb{Z} : m+j\frac{n-m}{2^i}<k\leq m+(j+1)\frac{n-m}{2^i}\right\}.
\end{equation*}
Then we introduce
\begin{equation*}
A_{i,j}=\prod_{k\in S_{i,j}} C_k,\quad P_{i,j}=\prod_{p\in S_{i,j}} p,
\end{equation*}
and
\begin{equation*}
R_{i,j}=M_m \prod_{0\leq r<j} A_{i,r} \pmod {P_{i,j}}.
\end{equation*}
For each prime $p\in S_{h,j}$ it follows that $r_p$ is congruent to the $(1,2)$-entry of the matrix $R_{h,j}$.

\section{Implementation}

For large integer arithmetic computations, we used the libraries GMP~\cite{Gr}
and NTT~\cite{Ha}. The NTT library supports multithreading without an
additional memory overhead and performs integer multiplication faster than GMP
routines when the operands are sufficiently large. This setup is similar to
the solution given in~\cite{CGH}. To generate a list of primes, we used the
implementation of the sieve of Eratosthenes provided by the FLINT
library~\cite{HJP}.
The source code of our implementation is available at
\url{https://github.com/milostatarevic/left-factorial}.

To compute all remainders $r_p$ for primes $p$ in an interval $(m,n]$,
we implemented the following four phases.

\subsection{Phase 1: computation of the $P_{i,j}$}

This phase consists of two parts. First, we generated a list of primes in $(m,n]$, then we
computed and stored all $P_{i,j}$ using a product tree. The time complexity of this
phase is $O((n-m+\sqrt{n})\log^{2+\varepsilon}n)$.

\subsection{Phase 2: computation of $M_m \pmod {P_{0,0}}$}

We computed $M_m$ by using a product tree. The time complexity of this phase
is $O(m \log^{3+\varepsilon}m)$. This phase represents the bottleneck of the
proposed algorithm.

In practice, each time we extended the computation to the next interval $(m,n]$, we
reused the intermediate multiplication results we stored from the previous iteration.
This approach allows us to reduce the computation time by a constant factor.
The optimal results were achieved when the stored values were just slightly less
than $2P_{0,0}$, which additionally required that the tree be partitioned carefully.

Unfortunately, this approach significantly increased space requirements
(measured in terabytes). As the data storage solutions are less expensive
compared to the cost of RAM or the price per CPU core, we decided to reduce
the computation time on account of the increase of the storage space. To
reduce the hard disk I/O we used a smaller solid-state drive, where we stored
intermediate results. This way the disk I/O did not represent a bottleneck.

As this phase is the most time expensive, the best performance is obtained if
the interval $(m,n]$ is as large as possible, which is limited by the
available RAM.

\subsection{Phase 3: computation of the $A_{i,j}$}

To compute $A_{i,j}$, we also used a product tree. To optimize space usage, we
stored only $A_{i,j}\pmod {P_{i,j}}$. Additionally, we used the results from
this phase to prepare the computation of $M_m$ for the next search interval as
described in Phase 2. The time complexity of this phase is
$O((n-m)\log^{3+\varepsilon}n)$.

\subsection{Phase 4: computation of the $R_{i,j}$}

This phase is similar to Phase 3, with the difference that we performed the
computation starting from the top level of the product tree $i=0$, going down
to the level $i=h$. The only values we had to store during this process
belonged to the level we were currently processing and those contained in one
level above. The time complexity of this phase is also
$O((n-m)\log^{3+\varepsilon}n)$.

\subsection{Verification of the results}

To verify a subset of computed values $r_p$, we used a procedure based on the
algorithm described in \cite{BGS}, with the time complexity
$O(p^{0.5+\varepsilon})$ per prime~$p$. The polynomial multiplication is
performed by using the NTL library \cite{Sh}.

\subsection{Hardware}

The computation was performed using a 6-core CPU (i7 6800K). The configuration
was equipped with 64GB of RAM and 16TB of disk space. The entire computation
took approximately $33\, 000$ core hours, where about $65$\% of the time was
spent in Phase 2. For a couple of the last blocks we processed, the time spent
in Phase 2 was approaching $80$\%.

\section{Results}

We calculated and stored $r_p$ for all primes $p$ less than $2^{40}$.
Heuristic considerations suggest that $!p$ is a random number modulo $p$ with
uniform distribution, so the probability that $r_p$ has any particular value
is approximately $1/p$, and the sum of reciprocals of the primes diverges.
Thus, we might expect that the probability to find a counterexample in an
interval $(2^m,2^n)$ is approximately $1-m/n$, and the expected number of
primes $p$ with $\abs{r_p}<\ell$ is approximately $(2\ell-1)\log(n/m)$
\cite{AT1}.

The new search covered the interval $(2^{34},2^{40})$. 
The above heuristic predicts that our chances of finding a counterexample were approximately~15\%. Although we
only found 24 primes with $\ell=100$ in our interval in comparison with the
expected value 32, the heuristics give good estimates for higher values
of~$\ell$. For $\ell=10\,000$ the expected number of primes in this interval
is $3250$, which is close to the actual value~$3237$. Similarly, for
$\ell=10\,000\,000$, the predicted value~$3\,250\,379$ is close to the actual
value~$3\,250\,456$. The results for $\abs{r_p}<100$ are presented in ~Table~\ref{Table1}.

Additionally, we have used our new algorithms in a search for \emph{socialist
primes}. The socialist primes are the primes~$p$ for which the residues of
$2!$, $3!$, \dots, $(p-1)!$ modulo $p$ are all distinct \cite{Tr}. Erd\H{o}s
conjectured that no prime $p>5$ is socialist, see ~\cite[Section~F11]{Gu}.
In our previous work \cite{AT2} we proved that a socialist prime $p$ needs to
satisfy $(!p-2)^2\equiv 1\pmod p$, and showed there are no socialist primes
with $5<p<10^{11}$. Our new results confirm that there are no primes $p$ in
the interval $(2^{34},2^{40})$ such that the remainder $r_p$ satisfies the
desired congruence. Consequently, there are no socialist primes $p$ with
$5<p<2^{40}$.  

\begin{center}
\begin{table}[t]
\begin{tabular}{||c|c||c|c||c|c||}
\hline
$p$ & $r_{p}$ & $p$ & $r_{p}$ & $p$ & $r_{p}$\\
\hline
$22\,370\,028\,691$ & $-55$ & $153\,736\,627\,747$ & $24$ & $450\,798\,203\,041$ & $52$ \\
\hline
$34\,212\,035\,633$ & $47$ & $203\,109\,046\,969$ & $-73$ & $541\,389\,125\,113$ & $-9$ \\
\hline
$35\,420\,262\,113$ & $-24$ & $252\,164\,235\,031$ & $84$ & $576\,365\,852\,729$ & $5$ \\
\hline
$39\,541\,338\,091$ & $-1$ & $296\,599\,719\,739$ & $-67$ & $581\,743\,725\,197$ & $28$ \\
\hline
$71\,848\,806\,989$ & $-87$ & $315\,631\,019\,399$ & $72$ & $668\,487\,297\,869$ & $-92$ \\
\hline
$94\,844\,067\,751$ & $-59$ &$342\,077\,311\,241$ & $-85$ & $740\,405\,032\,753$ & $-24$ \\
\hline
$102\,281\,886\,901$ & $19$ & $348\,036\,477\,379$ & $-77$ & $817\,880\,148\,803$ & $-46$ \\
\hline
$141\,853\,427\,273$ & $95$ & $425\,430\,768\,359$ & $9$ & $885\,831\,128\,921$ & $-35$\\
\hline          
\end{tabular}  
\vspace{0.2cm}
\caption{{The 24 primes in the interval $(2^{34},2^{40})$ such that  $\abs{r_p}<100$.}}    
\label{Table1}
\end{table}
\end{center}
\vspace{3mm}   

\section{Remarks}

After our article appeared in preprint, we learned that the theoretical aspects of
using remainder trees to compute the left factorial residues were also covered
in Rajkumar's master's thesis \cite{Ra}. We encourage the readers to read it.
Let us note that our work is independent and was published online at
approximately the same time as Rajkumar's work.
The computations and the results we presented in our paper are going back to 2017
and were initially presented at 14th Serbian Mathematical Congress in 2018 \cite{An}.         

\medskip {\bf Acknowledgements.} This work was partially supported by the Serbian Ministry of Education and Science, project No. 174012.
We would like to thank the referees for their careful reading and useful suggestions.

\bibliographystyle{amsplain}

\end{document}